\documentclass[12pt]{article}
\usepackage{amsthm,amsmath,amssymb,latexsym}
\title{Complete subvarieties in moduli spaces of rank 2 stable sheaves
on smooth projective curves and surfaces}
\author{Cristian Anghel}
\date{06. 03. 2001}
        \newtheorem{lem}{Theorem}[section]
        \newtheorem{thl}[lem]{Theorem}
        \newtheorem{cor} {Corollary}[section]
        \newtheorem{th2}[cor]{Theorem}
\begin{document}

\maketitle
\noindent {\small
            Institute of Mathematics of the Romanian Academy,
            P.O. Box 1-764, RO-70700 Bucharest, Romania}

        \vspace{10mm}

\section*{} Mathematics Subject Classification $(\ 2000\ )\ 14D20\ ,14J60$
\section*{Running title} Complete subvarieties in moduli spaces
        \vspace{10mm}

\noindent {\small {\bf Abstract} The aim of this paper is to prove
the existence of large complete subvarieties in moduli spaces of
rank two stable sheaves with arbitrary $c_{1}$ and sufficiently
large $c_{2}$, on algebraic surfaces. Then we study the
restriction of these sheaves  to curves of high
degree embedded in the surface. In the final section we gives a
relation with the spin strata defined by Pidstrigach and Tyurin.}

\section{Introduction}

Moduli spaces of rank 2 stable vector bundles over algebraic
surfaces had atracted considerable interest in the last years in
relation with the powerfull Donaldson's invariants. Apart their
detalied structure, a natural problem is to find a "reasonable"
bound for $c_{2}$ from which the moduli spaces are nonempty. The
answer of this question has been given by many authors in
different degrees of generality in \cite{A},\cite{H-L},\cite{Q}.

A second question related more closely to the Donaldson's
invariants is to find, at last in some special cases, complete
subvarieties of moduli spaces, which can in principle be used for
the calculus of Donaldson intersection numbers. For example in
\cite{B}, \cite{H} the authors construct complete curves in moduli
spaces of stable bundles with rank 2 and $c_{2}$ sufficently
great.

The aim of this paper is to present a construction of complete
subvarieties with large dimension (of order $2c_{2}$), in moduli
spaces of rank 2 stables sheaves with $c_{2}$ sufficently great,
extending the basic construction of \cite{A}.

Then we compare the dimension of our subvarieties with the O'Grady bound from
\cite{O'G} and
 we study the restriction of the
sheaves in our subvarieties to curves of high degree embedded in
the surface, obtaining also complete subvarieties in moduli spaces
of stables rank 2 bundles over curves. In the final section we
give the relation of ours subvarieties with the spin strata
defined in \cite{T}, \cite{P-T} for the spin-polinomial
invariants.

\section{Construction of subvarieties}

For fixing the notations let $S$ be a smooth projective surface, $H$ a very
ample polarisation and $L \in \mbox{Pic}(S)$ a line bundle on $S$.
Let $\widehat M_H(2,L,c_{2})$ be the moduli space of $H$-stable rank 2
torsion free sheaves $E$ with $det(E)=L$ and $c_{2}(E)=c_{2}$.
Let $n_{L}$ and $\beta(L,H)$ be defined as follows: for $L \cdot H > 0$
\[
        \begin{array}{ccl}
  n_{L}      & = &  \mbox{inf}\;  \{ n \in {\bf N}^{\ast} \mid
                        nH^{2} > \frac{L \cdot H}{2}, \;
                        nH^{2} > (L + K_{S}) \cdot H \}  \\

\alpha (L,H) & = & \mbox{max}\;
                \{h^{0}({\cal O}_{S}(L+K_{S}))+1, 2 + n_{L}H \cdot (L+K_{S}),
                1+n_{L} \frac{L\cdot H}{2} \} \; ,
        \end{array}
\]
and for arbritrary $L \cdot H$ :
\[
\beta (L,H) =  \left\{
  \begin{array}{lcc}
          \alpha (L,H)            &      \mbox{if}  &   L \cdot H > 0 \\
          \alpha (-L,H)           &      \mbox{if}  &   L \cdot H < 0 \\
          \alpha (L+2H,H) - H^{2} &      \mbox{if}  &   L \cdot H = 0
  \end{array} \right.
\]

In \cite{A} it was proven that the following theorem holds:

\begin{lem}
For all $c \geq \beta(L,H)$ there existe a rank two vector bundle
$H$-stable with $c_{1}(E) = L$ and $c_{2}(E) = c$.
\end{lem}

In the present note we are concerned with the following result:

\begin{thl}
For $c_{2}$ sufficently great (explicitely given) , the moduli space
$\widehat M_H(2,L,c_{2})$ contains a smooth complete subvariety
of dimension
\[
2c_{2}+h^0(-L)-\chi(-L)-1
\]
which has the structure of a projective fiber bundles over symmetric
product of a curve in $S$.
\end{thl}

The expression "sufficently great" means that $c_{2}$ must satisfy
the following conditions:
\[
c_{2} \geq \beta(L,H)
\]
\[
c_{2} > n_{L} \cdot L \cdot H
\]
\[
c_{2} \ge 2g(C)-1
\]
where $C \in \mid n_{L}H \mid $ is smooth.

{\it Proof:}

As in the proof of the {\bf Theorem 2.1.} in  \cite{A}, the principal
case is   $L \cdot H > 0$; the other situations are reduced to
this by dualising if  $L \cdot H < 0$ or by taking the tensor
product with ${\cal O}_{S}(H)$ if $L \cdot H = 0$. For simplicity
we assume in the sequel that $L \cdot H > 0$.

First of all, for any smooth curve $C$ on $S$ and any $Z \in
Div^{c_{2}}(C)$,
$Z$ can be seen in a single way as a $0$-dimensional subscheme in $C$ and
through the embedding $C\subset S$ as a $0$-dimensional local complete intersection
subscheme on $S$.
So we have a closed embedding of $Div^{c_{2}}(C)= Sym^{c_{2}}(C)$ in
$Hilb^{c_{2}}(S)$-the Hilbert scheme of $0$-dimensional subscheme of
length $c_{2}$ in $S$.

The main idea is that for a smooth curve $C \in \mid n_{L}H \mid $
and arbitrary $Z \in Div^{c_{2}}(C)$,any extension:
\[
0 \to {\cal O}_{S} \to E \to {\cal O}_{S}(L) \otimes {\cal J}_{Z} \to 0 \:
\]
is $H$-stable and the generic extension is locally free. The $H$-stability
is in fact a consequence of the basic lemma in \cite{A} and the
generic locally freenes follows from the general theory in  \cite{gr}.

Also, for $ c_{2} > n_{L} \cdot L \cdot H $ there are no two of
them isomorphic because the existence of an isomorphism between
the $E$'s of two different extention would imply that we have a
diagram of the following type:

\[
\begin{array}{ccccccccc}
      &       &       &       &
O
&       &       &       &       \\
      &       &       &       &
\uparrow
&       &       &       &       \\
      &       &       &       &
{\cal O}_S(L)J_{Z'}
&       &       &       &       \\
      &       &       &       &
\uparrow
&       &       &       &       \\
O  & \longrightarrow
   & {\cal O}_S
   & \longrightarrow
   & E
   & \longrightarrow
   & {\cal O}_S(L)J_{Z}
   & \longrightarrow
   & 0                          \\
      &       &       &       &
\uparrow
      &
\nearrow  \psi
      &       &       &         \\
      &       &       &       &
{\cal O}_S
&       &       &       &       \\
      &       &       &       &
\uparrow
&       &       &       &       \\
      &       &       &       &
0
&       &       &       &
\end{array}
\]

The arrow $\psi$ dose not vanish because otherwise the two extentions would be the same.
This  show that  $L$ has a section vanishing on $Z$ and this contradict:

\[
c_{2} > n_{L} \cdot L \cdot H \  .
\]

These extensions are parametrised by a projective fiber bundle
over $Sym^{c_{2}}(C)$; the fiber over $Z$ is the projectivisation
of

\[
\begin{array}{c}
{\rm Ext}^1 (LJ_Z,\, {\cal O}_S) \cong H^1(S,K_SLJ_Z)^*
\end{array}
\]

The dimension of this bundle is easily calculeted using the cohomology sequence
of the standard sequence:

\[
0 \to {\cal O}_{S}(K+L) \otimes {\cal J}_{Z} \to {\cal O}_{S}(K+L) \to {\cal O}_{Z} \to 0 \:
\]

In cohomology we have the following:
\begin{eqnarray} \nonumber
0 \to H^0( {\cal O}_{S}(K+L) \otimes {\cal J}_{Z}) \to H^0({\cal O}_{S}(K+L)) \to
H^0({\cal O}_{Z})\to    \\
\to H^1( {\cal O}_{S}(K+L) \otimes {\cal J}_{Z}) \to H^1({\cal O}_{S}(K+L)) \to 0
\nonumber
\end{eqnarray}

So we have $dim H^1(S,K_SLJ_Z)= c_{2} +h^{1}(-L)-h^{2}(-L)=
c_{2}+h^{0}(-L)-\chi (-L)$ and consequently the projective bundle we have
constructed has the fiber of dimension

\[
c_{2}+h^{0}(-L)-\chi (-L)-1
\]
and the complete subvariety has the dimension announced in
{\bf Theorem 2.2.}

Furthermore, if $c_{2} \ge 2g(C)-1$ then every line bundle over
$C$ of degree $c_{2}$ is nonspecial and therefore $Sym^{c_{2}}(C)$
is itself a smooth projective fiber bundle over the jacobian
$J(C)$ of $C$, with the fiber of dimension $c_{2} - g(C)$.

Finaly, let us observe that if the line bundle $L$ is nef and big,
then the term $h^0(-L)$ vanish and so the dimension of our subvarieties
depends only on numerical variables.

\section{Comparison with the O'Grady bound}

As we have seen in the preceeding section, the subvarieties we have
constructed contain sheaves which are not locally free. The aim
of this section is to relate the dimension of our subvarieties
with a result obtained by O'Grady in \cite{O'G} concerning the
maximal dimension of complete subvarieties in moduli spaces of
stables vector bundles.

For fixing the notations let $E$ a rank two sheaf over $S$. Let
\[
\Delta_{E}=c_{2}-c_{1}^{2}/4
\]
\[
\Delta_{0} (H) =  \left\{
  \begin{array}{lcc}

    3H^{2}                  &      \mbox{if}  &   K \cdot H < 0 \\
    3H^{2} (1+\frac{K \cdot H}{ H^{2}})^{2} &      \mbox{if}  &   K \cdot H \geq 0

  \end{array} \right.
\]
\[
\lambda_{2}=\frac{23}{6}
\]
\[
\lambda_{1}(H)=\frac{1}{2\sqrt{3H^{2}}}\cdot(4H^{2}+3K\cdot H+4)
\]
Finally let $ \lambda_{0}(H):= $
\[
\frac{3(K\cdot H+H^{2}+1)^{2}}{2H^{2}}+\frac{(K\cdot
H)^{2}}{4H^{2}}-\frac{K^{2}}{4}+4-3\chi({\cal O}_S), if  K \cdot H < 0
\]
and
\[
 \frac{3(K\cdot H)^{2}}{H^{2}} + 6 K \cdot H +\frac{3H^{2}}{2}-
\frac{K^{2}}{4} +8-3\chi({\cal O}_S), if  K \cdot H \geq 0 .
\]
Let $\mathcal{M}$ be the moduli space of rank two $H$-semistables
sheaves with $c_{2}$ and the determinant fixed. The theorem of
O'Grady we are concerned with, is the following:

\begin{thl}
If $H$ is very ample, $\Delta_{E} > \Delta_{0} (H)$ and
$\mathcal{V}$ is a complete subvariety of $\mathcal{M}$ such that:
\[
dim {\mathcal V}  >  \lambda_2 \cdot \Delta_{E} + \lambda_{1}(H) \cdot
\sqrt{\Delta_{E}}+ \lambda_{0}(H),
\]
then $\mathcal{V}$ contain sheaves which are not locally free.

\end{thl}

A very simple computation shows that, at least in the limit $c_{2}\gg
0$, the subvarieties constructed in {\bf Section 2}
have the dimension less than the bound in O'Grady theorem. So we
can expect that there are complete subvarieties of the same
dimension, as of ours, in which all the sheaves are locally free
and $H$-stables. Unfortunately we can not provide an explicit
construction for such ( potentially existing ) subvarieties.

\section{Application to moduli spaces over curves}

Let $D$ be a smooth projective curve of genus $g$ over the complex
field. The moduli space of rank two semi-stable vector bundles
with fixed determinant of degree $d$ is known to be a rational
projective variety (smooth if $d$ is odd or $g=2$) of dimension
$3g-3$. Its algebraic structure depends of $D$ but as a topological
space it is identified, if $d$ is even, with the space of $SU(2)$
representations of $\pi_{1}(C)$. Let denote this space by
${\mathcal{M}}_{g}(d)$, assuming in what follows that $d=0$ or
$1$.

We will apply the construction of {\bf Section 2} for obtaining for
certain values of $g$ a projective space $\mathbb P ^k$ embedded
in ${\mathcal{M}}_{g}(d)$.

For this, let's denote by $\mathcal{V}$ the complete variety
constructed in {\bf Section 2}. With the notations used there, let $Z$
be a reduced $0$-cycle on the curve $C \in \mid n_{L}H \mid $ formed
by $c_{2}$ distinct points. Let $l\gg 0$ an integer and $D\in \mid
lH \mid$ a smooth divisor which does not intersect $Z$. Let
${\mathcal{V}}_{Z}$ the fiber over $Z$ which is in fact $
{\mathbb{P}} ( {\rm Ext}^1 (LJ_Z,\, {\cal O}_S))$.

The main point is the following theorem which is a version of the
{\bf Theorem 1.1.} in \cite{tyurin} adapted for the presence of non
locally free sheaves.

\begin{thl}
For $l\gg 0$ and  $D\in \mid lH \mid$ as above, the restriction of
every sheaf $E\in {\mathcal{V}}_{Z}$ on $D$ is stable and the
morphism
\[
{\mathcal{V}}_{Z} \rightarrow {\mathcal{M}}_{D}(d)
\]
is an embedding, where $d=L\cdot D$.
\end{thl}

We will give the sketch of the proof for the case we are interested in, where
there are non locally free sheaves.

{\it Proof:}Let $ M_H(2,L,c_{2})$ be the moduli space of rank two
$H$-stable vector bundles with the prescribed determinant and
second Chern class. The result of Tyurin cited above asserts that
for any $c_{2}$ there exists a $l_{0}(c_{2})$ such that for any
$l\geq l_{0}(c_{2})$ and a generic smooth $D\in \mid lH\mid$ the
restriction map
\[
\bigcup M_H(2,L,k) \rightarrow {\mathcal{M}}_{D}(d)
\]
is an embedding, where the union is taken for all $k\leq c_{2}$
and $d=L\cdot D$.

In what follows we will choose a $l$ and a $D$ as in the Tyurin
theorem and such that $D \cap Z= \emptyset$. Therefore if $E\in
\mathcal V _{Z}$ is a locally free sheaf, the restriction $E_{|D}$
is stable and, if ${\mathcal{V}}_{Z,lf}$ is the locus corresponding
to locally free sheaves, then
\[
{\mathcal{V}}_{Z,lf} \rightarrow {\mathcal{M}}_{D}(d)
\]
is an embedding, by Tyurin theorem. We are therefore interested
for the behaviour of the restriction map on the non locally free
locus. Let $E$ be such a sheaf. As we are on a surface, there is
an exact sequence of the following type:
\[
0 \to E \to E^{\vee \vee} \to Ct(E) \to 0 \:
\]
where $ E^{\vee \vee}$ is the bidual of $E$ and $Ct(E)$ is the
cotorsion sheaf. It is well known that the followings facts are
true: $E^{\vee \vee}$ is locally free, $H$-stable if and only if
$E$ is $H$-stable,
 $c_{2}(E^{\vee \vee}) = c_{2}(E) - lenght(Ct(E))$, $Ct(E)$ is a
sheaf with a $0$-dimensional support included in the locus where
$E$ is not locally free therefore in $Z$, and $E$ is isomorphic
with $ E^{\vee\vee}$ on $S\backslash Z$. So, as $D \cap Z=
\emptyset$ we conclude that
\[
 E_{|D}\simeq  E^{\vee \vee}_{|D}
\]

But $ E^{\vee \vee}$ is $H$-stable, locally free and from above
$c_{2}(E^{\vee \vee}) < c_{2}(E)$. So by Tyurin theorem $ E^{\vee
\vee}_{|D}$ is stable. Moreover, for any $E_{1}\in
{\mathcal{V}}_{Z}$ and $E_{1} \neq E$ we have ${E_1}_{|D} \neq
E_{|D} $ for the following reason: if $E_{1}$ is locally free we
apply Tyurin theorem to $E_{1}$ and $ E^{\vee \vee}$; if $E_{1}$
is not locally free we apply again Tyurin theorem to $E_{1}^{\vee
\vee}$ and $ E^{\vee \vee}$.QED

$\mathbf{Remarque}$.The presence of such linear spaces in the moduli spaces
of stables bundles over curves is not surprising in view of the rationality
of these spaces proved in \cite{sch}.

\section{Relation with the spin strata}

We recall the definition of spin strata from \cite{T}, \cite{P-T}
in the algebraic geometric context, without considering the
differential geometric counterpart.

Let $\mathcal{M}$ be the moduli space of rank two $H$-semistables
sheaves with $c_{2}$ and determinant $L$ fixed. The spin strata
${\mathcal{M}}_{i}$ are defined as follows:
$$
{\mathcal{M}}_{i}=\{E\in {\mathcal{M}} | H^{0}(E)+ H^{2}(E) \geq
i\} $$.
Moreover, the strata ${\mathcal{M}}_{i}$ has a
decomposition
$$ {\mathcal{M}}_{i}= \bigcup {\mathcal{M}}_{j,l} \ ,$$
where the union is taken for all $j,l$ such that $j+l=i$ and the
${\mathcal{M}}_{j,l}$ are defined by
$$ {\mathcal{M}}_{j,l}=\{E\in
{\mathcal{M}} | H^{0}(E) \geq j and H^{2}(E) \geq l\}
$$.

Using the construction in {\bf Section 2 } in the case $L \cdot H > 0$ we
obtain the following corollary of {\bf Theorem 2.2}:
\begin{cor}
If $L \cdot H > 0$ the spin stratum ${\mathcal{M}}_{1}$ contains a smooth
complete subvariety $\mathcal{V}$ of dimension
$2c_{2}+h^0(-L)-\chi(-L)-1$.

Moreover, if $H^{2}({\mathcal{O}}_{S})=H^{2}(L)=0$, $\mathcal{V}$ is not
included in any higher codimensional stratum ${\mathcal{M}}_{i}$ for
$i\geq 2$ and in fact it is contined in ${\mathcal{M}}_{1,0}$.

\end{cor}

{\bf Proof}: The first assertion is obvious from the construction
of the bundles $E \in \mathcal{V}$.

For the second we must show that $H^{0}(E)=1$ and $H^{2}(E)=0$
in the hypothesis of the Corollary. Using the sequence
\[
0 \to {\cal O}_{S} \to E \to {\cal O}_{S}(L) \otimes {\cal J}_{Z} \to 0 \:
\]
we obtain that $H^{2}(E)$ is part of the following sequence:
\[
\to H^{2}({\cal O}_{S}) \to H^{2}(E) \to H^{2}({\cal O}_{S}(L) \otimes {\cal
J}_{Z}).
\]
Now, the first term is $0$ while the third is part of the sequence:

\[
0 \to H^{2}({\cal O}_{S}(L) \otimes {\cal J}_{Z}) \to H^{2}(L) \to 0
\]
obtained by taking the cohomology of the following:
\[
0 \to {\cal O}_{S}(L) \otimes {\cal J}_{Z} \to L \to {\cal O}_{Z} \to 0
\]
So, using that $H^{2}(L)=0$ we obtain $H^{2}(E)=0$ as we wished.

Suppose now that $H^{0}(E)\geq 2$. We obtain a diagram of the
following type:

\[
\begin{array}{ccccccccc}

O  & \longrightarrow
   & {\cal O}_S
   & \longrightarrow
   & E
   & \longrightarrow
   & {\cal O}_S(L)J_{Z}
   & \longrightarrow
   & 0                          \\
      &       &       &       &
\uparrow
      &
\nearrow  \psi
      &       &       &         \\
      &       &       &       &
{\cal O}_S
&       &       &       &       \\
      &       &       &       &
\uparrow
&       &       &       &       \\
      &       &       &       &
0
&       &       &       &
\end{array}
\]
and so $\psi$ is a nonzero section in ${\cal O}_{S}(L) \otimes
{\cal J}_{Z}$. But this contradicts one of the fundamental
hypothesis of our construction concerning the choice of $c_{2}$
and of $Z$,namely:
\[
c_{2} > n_{L} \cdot L \cdot H.
\]
So the Corollary is proved. QED

In the rest of this section we are concerned with the construction
of complete subvarieties which are not strictly contained in any
spin strata ${\mathcal{M}}_{j,l}$ with $j\geq 1$ or $l\geq 1$, and
so we must modify the construction in {\bf Section 2} for
obtaining bundles without sections. For this let's making the
following notation:
\[
L'=L+2mH
\]
\[
{c'}_{2}=c_{2}+mH\cdot L +m^{2}H^{2}
\]
where $m$ is a positive integer. We want to apply the main
construction in $\mathbf{Section2}$ for the modified Chern classes
$L'$ and ${c'}_{2}$. So we will obtain stables sheaves $E'$ with
the above Chern classes which are setting in exact sequences of
the following type:
\[
0 \to {\cal O}_{S} \to E' \to {\cal O}_{S}(L+2mH) \otimes {\cal J}_{Z'} \to 0 \:
\]
where $Z'$ is a zero-dimensional subscheme of length $l(Z')= {c'}_{2}$ such
that the conditions in {\bf Theorem 2.2.} are satisfied.

The bundles $E=E'\otimes {\cal O}_{S}(- mH)$ are obviously
stable, as the $E'$ are, and with Chern classes $L$ and $c_{2}$.
So we obtain a new family ${\mathcal{V}'}_{m}$ of stable bundles
with the prescribed Chern classes.

The above remarks can be summarized in the following Corollary:
\begin{cor}

For $c_{2}\gg 0$ the moduli space $\mathcal{M}$ contains the complete
subvariety ${\mathcal{V}'}_{m}$, the member of which are setting in
sequences of the following type:
\[
0 \to {\cal O}_{S}(-mH) \to E \to {\cal O}_{S}(L+mH) \otimes {\cal J}_{Z'} \to 0 \:
\]

\end{cor}

{\bf Remark.} The explicit meaning  of $c_{2}\gg 0$ and the
dimension of ${\mathcal{V}'}_{m}$ can be obtained as in {\bf Section 2} by
making the computation for $L'$ and ${c'}_{2}$ and then translating in terms
of $L$ and $c_{2}$.

Another significant fact is that for fixed $c_{2}\gg 0$ there are only a
finite number of values for $m$ such that the ${\mathcal{V}'}_{m}$ can be
constructed by the above procedure. This is a consequence of the fact that
for $m\gg 0$ the growth of ${c'}_{2}$ is as $m^{2}H^{2}$, while the growth
of $\beta (L',H)$ is at last as $2m^{2}H^{2}$ and the fundamental
condition
\[
{c'}_{2} \geq \beta(L',H)
\]
forbid that $m$ would be too large.

Also, concerning the dimension of ${\mathcal{V}'}_{m}$, a simple computation
shows that it is of the same order as the dimension of $\mathcal{V}$, that
is of order $2c_{2}$. The topology of ${\mathcal{V}'}_{m}$ is a little changed
in that now it is a projective bundle over a symmetric power of
order $c_{2}+mH\cdot L +m^{2}H^{2}$ over a curve in $S$.$\Box$

We want now to see the relation of the ${\mathcal{V}'}_{m}$ with
the spin strata. First of all if $E\in {\mathcal{V}'}_{m}$ we have
$H^{0}(E)=0$. If not, we have a diagram:

\[
\begin{array}{ccccccccc}

O  & \longrightarrow
   & {\cal O}_{S}(-mH)
   & \longrightarrow
   & E
   & \longrightarrow
   & {\cal O}_{S}(L+mH) \otimes {\cal J}_{Z'}
   & \longrightarrow
   & 0                          \\
      &       &       &       &
\uparrow
      &
\nearrow  \psi
      &       &       &         \\
      &       &       &       &
{\cal O}_S
&       &       &       &       \\
      &       &       &       &
\uparrow
&       &       &       &       \\
      &       &       &       &
0
&       &       &       &
\end{array}
\]
with $\psi$ obviously nonzero because ${\cal O}_{S}(-mH)$ has no
sections. So $\psi$ gives a nonzero section in ${\cal O}_{S}(L+mH)
\otimes {\cal J}_{Z'}$ which is impossible by the choice of $Z'$
and by the fundamental condition

\[
{c'}_{2} > n_{L'} \cdot L' \cdot H
\]
used in the construction of the $E'$.

With the precedings  notations we are now able to state
the main result of this section:

\begin{th2}

If $L\cdot H>0$ and $L\cdot H>K\cdot H$ then
the complete subvariety ${\mathcal{V}'}_{m}$ in $\mathcal{M}$ is not
strictly contained in any spin stratum
${\mathcal{M}}_{j,l}$ with $j\geq 1$ or $l\geq 1$. More exactly for
$E \in {\mathcal{V}'}_{m}$ generic we have:
\[
H^{0}(E)=H^{2}(E)=0
\]
\end{th2}

$\mathbf{Proof}$: The conclusion for the $H^{0}$ follows for all the $E$
in ${\mathcal{V}'}_{m}$ by the preceedings considerations.

In what follows consider $E$ in ${\mathcal{V}'}_{m}$ generic so that it
is locally free and the zero-dimensional sub-scheme $Z'$ which
correspond to it consists of distinct reduced points. By taking the
dual of the defining sequence of $E$ and by tensoring it with
${\cal O}_{S}(K)$ we obtain the following:

\[
0 \to {\cal O}_{S}(-L-mH+K) \to E^{\vee} \otimes {\cal O}_{S}(K)
\to {\cal O}_{S}(K+mH) \otimes {\cal J}_{Z'} \to 0 \:
\]
By Serre duality we must  show that $H^{0}(E^{\vee} \otimes
{\cal O}_{S}(K))=0$. Suppose the contrary. We obtain the following
diagram:
{\footnotesize
\[
\begin{array}{ccccccccc}

O  & \longrightarrow
   & {\cal O}_{S}(-L-mH+K)
   & \longrightarrow
   & E^{\vee} \otimes {\cal O}_{S}(K)
   & \longrightarrow
   & {\cal O}_{S}(K+mH) \otimes {\cal J}_{Z'}
   & \longrightarrow
   & 0                          \\
      &       &       &       &
\uparrow  \varphi
      &
\nearrow  \psi
      &       &       &         \\
      &       &       &       &
{\cal O}_S
&       &       &       &       \\
      &       &       &       &
\uparrow
&       &       &       &       \\
      &       &       &       &
0
&       &       &       &
\end{array}
\]
}
Now, if $\psi =0$ we obtain a section in ${\cal O}_{S}(-L-mH+K)$
which is impossible because the hypothesis $L\cdot H>K\cdot H$
would imply
\[
(-L-mH+K)\cdot H < 0
\]
So $\psi$ is not zero and it defines a section in ${\cal O}_{S}(K+mH) \otimes {\cal
J}_{Z'}$. But in the fundamental construction we have
\[
{c'}_{2}=lenght(Z') > n_{L'} \cdot L' \cdot H
\]
and so $Z'$ can not lying on a curve in the linear system $|mH+K|$
because the hypothesis $L\cdot H>K\cdot H$ imply
\[
(mH+K)\cdot (n_{L'}H) < (mH+L)\cdot (n_{L'}H) < (L')\cdot (n_{L'}H)
\]

So for $E$ generic in ${\mathcal{V}'}_{m}$ we have
$H^{2}(E)=0$ and the theorem is proved. QED

\section*{}  Acknowledgements. The work to this paper was partially  supported by
a DFG grant in the Oldenburg University. I would like to thank
Prof. U. Vetter and Prof. N. Manolache for their interest, for the
stimulating atmosphere and for their very kind hospitality.

\end{document}